\documentstyle{amsart}

\def\del{\partial}

\def\ini{\operatorname{in}}

\def\into{\hookrightarrow}

\def\N{{\Bbb N}}
\def\sbe{\subseteq}

\def\var{\operatorname{V}}
\def\C{{\mathbb C}}
\def\D{{\mathcal D}}

\def\O{{\mathcal O}}

\newcounter{plm}
\newtheorem{thm}[plm]{Theorem}

\theoremstyle{definition}
\newtheorem{alg}[plm]{Algorithm}
\newtheorem{rem}[plm]{Remark}
\newtheorem{ex}[plm]{Example}
\newtheorem{notn}[plm]{Notation}
\def\bar#1{\overline{#1}}
\def\action{\bullet}

\begin{document}

\title{A Localization Algorithm for $D$-modules} 
\author[Toshinori Oaku, Nobuki Takayama and Uli Walther]
 {Toshinori Oaku$^{1}$, Nobuki
Takayama$^{2}$ and Uli Walther$^{3}$}
\address{$^1$Yokohama City University}
\email{oaku@@math.yokohama-cu.ac.jp}
\address{$^2$Kobe University and MSRI}
\email{takayama@@math.kobe-u.ac.jp}
\address{$^3$University of Minnesota and MSRI}
\email{walther@@math.umn.edu}
\date{November 4, 1998}
\begin{abstract}
We present a method to compute the holonomic extension of a $D$-module
from a Zariski open set in affine space to the whole space.
A particular application is the localization of coherent $D$-modules
which are holonomic on the complement of an affine variety.
\end{abstract}
\maketitle

Throughout this article let
$R=K[x_1,\ldots,x_n]$ be the ring of polynomials in $n$ variables
over the field $K$ of characteristic zero and 
$D=R\langle \del_1,\ldots,\del_n\rangle$ the $n$-th Weyl
algebra over $K$. Here, $\del_i=\frac{\del}{\del x_i}$. Suppose $f\in
 R$, let $J$ be  a left ideal of $D$ and set $M=D/J$. Let $X=K^n$,
$Y=\var(f)\sbe X$, $U=X\setminus Y$ and $j:U\into X$ the natural
inclusion.

In the landmark paper \cite{K2} M.\ Kashiwara proved that $\int_jj^{-1}M$ is
holonomic provided that  $M$ is
holonomic on $U$, but not necessarily on all of $X$.
In this note we make this important result algorithmic in the sense
that we provide an algorithm that computes the module
structure of $\int_jj^{-1}M$ over $D$.

We remark that $U$
is affine and $X=K^n$, and so $\int_jj^{-1}M=R_f\otimes_RM$ as
$D$-module where $R_f=R[f^{-1}]=\Gamma(\O_X,U)$.
An algorithm to compute $R_f\otimes_RM$ as a $D$-module has been given
in \cite{Oa3} under the assumption that $M$ is holonomic on $K^n$
and $f$-torsion free, and in \cite{O-T2} under the assumption that
$M$ is holonomic on $K^n$. The advantage of our new algorithm,
besides removing the requirement of holonomicity on $\var(f)$, is that
the natural map $M\to R_f\otimes_R M$ can be traced, which is an important
property for other computations as local cohomology
(\cite{Oa3}, \cite{W1}) and de Rham cohomology (\cite{O-T1}, \cite{W2}).

The algorithm is illustrated in a section devoted to computational
examples. 

\section{Description of the Algorithm}
We shall start with a motivating construction which will also provide
a skeleton of the proof of the correctness of our
algorithm. 

Assume that $K$ is in fact the field $\C$ of complex numbers.
Let $s, v$ be two new variables and let  $D_v=D\langle v,\del_v\rangle$.
\begin{notn}
Throughout we will denote $\partial f/\partial x_i$ 
by $f_i$.
If $P\in D_v$ is an operator interpreted as a polynomial in
$\del_1,\ldots,\del_n$, we will simply write $P(\del_x)$.
Similarly we will abbreviate
$P(\del_1-v^2f_1\del_v,\ldots,\del_n-v^2f_n\del_v)$ by
$P(\del_x-v^2f_x\del_v)$. 
\end{notn}
Consider the ring $D_v$ and the $\C$-vectorspace
$R_f[s]\otimes f^{-s}$. Here, $f^{-s}$ should be
thought of as a symbol that behaves like the complex valued function
$f^{-s}$ under differentiation and the tensor is over $\C$. 
We shall turn this vectorspace into a
module over $D_v$ by means of the action $\action$ defined as
follows. 
We require that all $x_i$ act as
multiplication on the left factor, 
while $\del_i$ acts by the ``product rule''
$\del_i\action(\frac{g(x,s)}{f^m}\otimes
f^{-s})=\del_i(\frac{g(x,s)}{f^m})\otimes
f^{-s}+\frac{g(x,s)}{f^{m+1}}(-s)f_i\otimes f^{-s}$. 
Thus, we only need to define the
action of $v$ and $\del_v$.
We set
\begin{eqnarray*}
v\action\left(\frac{g(x,s)}{f^m}\otimes
f^{-s}\right)&=&\frac{g(x,s+1)}{f^{m+1}}\otimes f^{-s},\\
\del_v\action\left(\frac{g(x,s)}{f^m}\otimes
f^{-s}\right)&=&-(s-2)\frac{f\cdot g(x,s-1)}{f^{m}}\otimes f^{-s}.
\end{eqnarray*}
We observe that $v$ is simply shifting $s$ up by one (in particular,
the extra $f$ in the denominator comes from shifting the exponent in
$f^{-s}$) and $\del_v$ is shifting $s$ down by one, together with the
``differential'' $2-s$ of the shift operator
(cf.\ the action defined in \cite{M}).

Armed with this action on $R_f[s]\otimes f^{-s}$ we can define
an action of $D_v$ on $D/J\otimes_R (R_f[s]\otimes f^{-s})$ as
follows. $v$ and $\del_v$ commute with the first factor in the
tensor product, $x_i$ acts by left multiplication on the first factor
and $\del_i$ acts by the product rule.

We detect a certain set of trivial relations for this
action. Namely, let us ask which operators on $D_v$ are 
candidates for annihilating the element $\bar 1\otimes 1\otimes
f^{-s}$, the bar denoting cosets in $D/J$. 
Independently of $J$, $1-vf$ is such an annihilator. Assume that 
$J$ is generated by $\{P_j(\del_x)\}_1^r$ where we think of $P_j$ as a
polynomial in $\del_x$ 
with left coefficients in $R$. 
Let us define a ring map $\phi$ from $D$ to $D_v$ by mapping $P\in D$ to
$\phi(P_j(\del_x))=Q_j(\del_x):=P_j(\del_x-v^2f_x\del_v)\in D_v$ (recall
$f_i=\del_i(f)$).  
It is easy to see
that this assignment preserves commutators and hence is an actual ring
map. We remark that $\phi$ represents the map
$U\stackrel{\cong}{\longrightarrow} \var(1-vf)\sbe \C^{n+1}$ given by
$x\to (1/f(x),x)$ on the level of differential operators, while the
inclusion $\phi(D)\into D_v$ corresponds to $\var(1-vf)\into \C^{n+1}$.

We claim that $Q_j(\del_x)$ annihilates $\bar 1\otimes 1\otimes
f^{-s}$. To prove this, we observe that $(\del_i-v^2f_i\del_v)\action
(\bar
P\otimes 1\otimes f^{-s})$ equals $\bar{\del_iP}\otimes 1\otimes
f^{-s}$ for all $P\in D$. Thus $Q_i(\del_x)\action(\bar 1\otimes 1\otimes
f^{-s})=\bar {P_i(\del_x)}\otimes 1\otimes f^{-s}=0$. We point out that
$-v\del_v+1$ acts as multiplication by $s$ and that this construction
is (formally) independent of the base field.

Our main statement is the following
\begin{thm}
Let $M=D/J$. 
Consider the left ideal $D_v\cdot\{1-xv,\phi(J)\}$ and the right ideal
$\{\del_v\}\cdot D_v$ in $D_v$. 
The quotient $D_v/(D_v\cdot\{1-xv,\phi(J)\}+\{\del_v\}\cdot D_v)$
is isomorphic to $M\otimes R_f$ as
$D$-module. Moreover, the natural map $M\to M\otimes_RR_f$ sends $\bar
1\in M$ to $\bar 1\otimes \frac{1}{f^2}$.
\end{thm}
It is worth noticing that this is in fact a
left $D$-module due to the fact that $D$ and $\del_v$ commute.

In the general context in which the theorem is stated, the $D$-module
$D_v/(D_v\cdot\{1-xv,\phi(J)\}+\{\del_v\}\cdot D_v)$ will not be a finitely generated
$D$-module. If however $M$ is holonomic on $U=K^n\setminus
\var(f)$, then 
$M\otimes_RR_f=D_v/(D_v\cdot\{1-xv,\phi(J)\}+\{\del_v\}\cdot D_v)$ will be
holonomic, and in 
particular finitely
generated. 
In order to compute the structure of the localized module,
we observe that it is the integration of the left $D_v$-module
$D_v/D_v\cdot\{1-xv,\phi(J)\}$ with respect to $\del_v$ 
(cf.\ \cite[Theorem 5.7]{Oa3}, \cite[Section 6]{O-T1}).

Methods for the algorithmic computation of $0$-th integration worked
out in \cite[Theorem 5.7]{Oa3}
(i.e., the restriction to $v=0$ of the Fourier transform of the module,)
together with our theorem give the following
algorithm. 
Let $N$ denote $D_v/(\{\del_v\}\cdot D_v +D_v\cdot\{1-fv,\phi(J)\})$.
\begin{alg}
\label{loc-alg}
Input: $f\in R$; $\{P_1,\ldots,P_r\}\subseteq D$ generating $J$, $D/J$
holonomic on $U=K^n\setminus \var(f)$.

\noindent Output: $k\in \N$; $\{Q_1,\ldots,Q_t\}\subseteq D$ such that
$D/D\cdot\{Q_1,\ldots,Q_t\}$ is isomorphic to $D/J\otimes _RR_f$ generated by
$\bar 1\otimes \frac{1}{f^{k+2}}$.

\noindent Begin.
\begin{enumerate}
\item For $P_i, i=1,\ldots,r$ 
compute $\phi(P_i)\in D_v$ defined by replacing $\del_x$ by
$\del_x-v^2f_x\del_v $.
\item Compute the $b$-function $b(s)$ for integration of
$D_v/D_v\cdot\{1-fv,\phi(J)\}$ with respect to $\del_v$. That is, find
$K[v\del_v]\cap\,\ini_{w}(D_v\cdot\{1-fv,\phi(J)\})$ where $w$ is the weight
assigning 1 to $v$, $-1$ to $\del_v$ and 0 to all other
variables. Replace $v\del_v$ by $-s-1$.
\item Let $k$ be the largest non-negative integer root of $b(s)$. 
If there is no such root, then output $I=D$, else continue.
The integral
$N$ is generated by the coset of
$v^k$ in $N$. 
\item Compute the annihilator $I$ over $D$ of $v^k$ in $N$.
\item $D/J\otimes_R R_f$ is generated by $\bar 1\otimes f^{-(k+2)}$ and
isomorphic to $D/I$.
\item Return $k$ and $I$.
\end{enumerate}
End.
\end{alg}
An algorithm to perform step 2 of the algorithm is given in
\cite[Algorithm 4.5]{Oa2}.

The steps 3 and 4 are nothing but an algorithm to get the $0$-th integral
of $D$-modules.
Here is a more precise description of steps  3 and 4.
\begin{enumerate}
\item Let $G = \{ g_1, \ldots, g_m \}$ be a Gr\"obner basis of
$D_v \cdot \{ 1-vf , \phi(J) \}$ with respect to $w$
(see, e.g.,  \cite[section 1.1]{SST}).
\item Let $G_k$ be the set
$$ \{\, {\rm normalForm}( v^i g_j , \{ \del_v \}\cdot D_v )\,|\,
    j=1, \ldots, m,  0 \leq i \leq k - {\rm ord}_w (g_j) \}.
$$
\item Regard $D\cdot \{G_k\}$ as a left submodule of the free module
$$ D\, v^0 + D\, v^1 + \cdots + D \, v^k = D^{k+1} $$
and find generators of
$D \cdot \{ G_k \} \cap D\, v^k$.
The last intersection can be computed by using an order to
eliminate $v^0, \ldots, v^{k-1}$,
i.e, by using an order $\succ$ such that
$ a_i v^i \succ a_k v^k$ for all $i=0, \ldots, k-1$ and 
$a_i, a_k \not= 0$ in $D$.
\end{enumerate}
Here, we put
$$ {\rm ord}_w (g):= \max\left\{ w \cdot (\alpha,\beta) \,|\, 
     g = \sum a_{ \alpha \beta } v^\alpha \del_v^\beta \right\}
$$
following the notation of \cite{SST}.
${\rm normalForm}(g,\{\del_v \}\cdot D_v)$ means taking the normal
form of $g$ with respect to the right ideal $\{\del_v\} \cdot D_v$.
For example,
$$ {\rm normalForm}(v^2 \del_v^2 ,\{\del_v \}\cdot D_v) =
   {\rm normalForm}(\del_v^2 v^2 - 4 \del_v v + 2 ,\{\del_v \}\cdot D_v) 
   = 2.
$$

\section{Proof of correctness of the algorithm}
Let us first provide the 
\par
\noindent{\bf Proof of the Theorem:} 
Given an operator $P\in D_v$ we shall define its {\em normal form} in $D_v$.
The goal is a presentation of $P$ where $\del_i-v^2f_i\del_v$ take the
position of $\del_i$. In other words, we are aiming for a sum of the
form $P=\sum \del_v^av^bp_{ab}(x)q_{ab}(\del_x-v^2f_x\del_v)$ where
$p_{ab} $ and $q_{ab}$ are polynomials in $n$ variables over $K$.
To this end, write first $P$ as $\sum
\del_v^av^bp_{ab}(x)q_{ab}(\del_x)$. 
The operator $P^1:=\sum\del_v^av^bp_{ab}(x)q_{ab}(\del_x-v^2f_x\del_v)$
has the property that $P-P^1$
will have lower degree in $\del_x$ than $P$. If $P-P^1=0$,
quit. Otherwise write $P-P^1$ as $\sum
\del_v^av^bp^1_{ab}(x)q^1_{ab}(\del_x)$ and set
$P^2:=\del_v^av^bp^1_{ab}(x)q^1_{ab}(\del_x-v^2f_x\del_v)$. Repeat
this procedure until we arrive at $P^l=0$. Then $P=P^1+\cdots+P^l$ and
this sum is the desired normal form in $D_v$.
We shall write $\tilde P$ for the normal form of
$P$ in $D_v$. 

The normal form in $D_v$ of an operator induces a normal form in 
$D_v/\{\del_v\}\cdot D_v$ by removing all terms in $\{\del_v\}\cdot D_v$ 
from the normal
form of $P$ in $D_v$. We denote the normal form of $P\in D_v$ in
$D_v/\{\del_v\} \cdot D_v$ by $\bar{\tilde P}$.
Of course, $P+\{\del_v\} \cdot D_v=\tilde P+ \{\del_v\} \cdot D_v$ as
cosets. We notice that both normal forms are unique.
As an example, consider the normal form of $\del_i$ in
$D_v/\{\del_v\}\cdot D_v$: 
$P=\del_i$, $P^1=\del_i-v^2f_i\del_v$, thus
$P^2=v^2f_i\del_v=\del_vv^2f_i-2f_iv$, $P^3=0$ and  hence
$\bar{\tilde P}=(\del_i-v^2f_i\del_v)-2f_iv+ \{\del_v\} \cdot D_v$. 

As a $D$-module, the quotient 
$N=D_v/(\{ \del_v\} \cdot D_v +D_v\cdot\{1-vf,\phi(J)\})$ is
spanned by the powers of $v$, which are in normal form. 
Let us try to understand the relations among these powers. To this
end note that
$(1-fv)(\del_i-v^2f_i\del_v)=(\del_i-v^2f_i\del_v+f_iv)(1-fv)$ which
implies that $\bar{(1-fv)\tilde P}=0$ in $N$ for all $P\in D_v$. 
This proves that $\bar{f(v\tilde P)}=\bar{\tilde P}=\bar P$ in
$N$ for all $P\in D_v$ and we conclude that every element in $N$
is divisible by $f$.

Moreover, the division by $f$ is unique since 
if $f\bar{\tilde Q}=f\bar{\tilde Q'}$ in $N$ then $\bar{vf(\tilde Q-\tilde
Q')}=0$ in $N$ and that 
implies $\bar{\tilde Q}=\bar{\tilde Q'}$ in $N$.
Conversely, for all $\bar{\tilde P}\in N$ there is a unique
$\bar{\tilde Q}\in N$ 
for which $\bar{v\tilde Q}=\bar{\tilde P}$ because then necessarily
$\bar{\tilde 
Q}=f\bar{\tilde P}$. So $f$ is an invertible operator on $N$ and the
inverse is given by $v$ acting on the normal form representative of a
coset. 

We shall now make use of the action we defined in the previous section.

The map $D_v\to D/J\otimes_R (R_f[s]\otimes f^{-s})$
sending $1$ to
$\bar 1\otimes  1\otimes f^{-s}$
defined via the action $\action$ is surjective. 
Consider the combined map $D_v\to
D/J\otimes_R( R_f[s]\otimes f^{-s})\to \left(D/J\otimes_R( R_f[s]\otimes
f^{-s})\right)/(s-2)\left(D/J\otimes_R( R_f[s]\otimes f^{-s}\right))$. 

This is a surjective map of left $D$-modules and the kernel contains
$\{\del_v\}\cdot D_v +D_v\cdot\{1-vf,\phi(J)\}$. If an operator in $D_v$
is in the kernel, we can assume it to be in normal form in 
$D_v/\{\del_v\}\cdot D_v$ 
since $\{\del_v\} \cdot D_v$
is in the kernel. Any such operator may be written as
$v^kP(\del_x-v^2f_x\del_v)$ 
where $P(\del_x)\in D$. We conclude that $P$ must already
be in the 
kernel since $v^kP(\del_x-v^2f_x\del_v)\action(\bar 1\otimes 1\otimes
f^{-s})=v^k\action(\bar{P(\del_x)}\otimes 1\otimes
f^{-s}=\bar{P(\del_x)}\otimes f^{-k}\otimes f^{-s}$. 
Thus $P(\del_x-v^2f_x\del_v)\action (\bar 1\otimes 1\otimes f^{-s})\in
(s-2)D/J\otimes_R(R_f[s] 
\otimes f^{-s})$ and therefore $f^jP(\del_x)\in D[s]\cdot
(J,(s-2))$ for some $j\in \N$. Since neither $P$ nor $J$ does 
contain $s$, setting $s=2$ we see that
$P(\del_x)$ is in the $f$-saturation of $J$, $f^jP(\del_i)\in J$.

$f^jP(\del_x)\in J$ implies that
$f^jP(\del_x-v^2f_x\del_v)\in\phi(J)$. Modulo $D_v\cdot\{1-vf\}$ this is equivalent 
to $P(\del_x-v^2f_x\del_v)\in\phi(J)$ since
$f^jP(\del_x-v^2f_x\del_v)$ is in normal form. This simply means
that $P$ is in the kernel of $D_v\to
D/J\otimes_R (R_f[s]\otimes f^{-s})\to \left(D/J\otimes_R( R_f[s]\otimes
f^{-s})\right)/(s-2)\left(D/J\otimes_R( R_f[s]\otimes f^{-s})\right)$ if
and only if it is in $\{\del_v\}\cdot D_v +D_v\cdot\{1-vf,\phi(J)\}$. 

Finally we observe that 
\[
\left(D/J\otimes_R( R_f[s]\otimes
f^{-s})\right)/(s-2)\left(D/J\otimes_R( R_f[s]\otimes f^{-s})\right)
\]
is
isomorphic to $D/J\otimes_R R_f$ by means of the map
$\bar P\otimes \frac{g(x,s)}{f^m}\otimes f^{-s}\to
\bar P\otimes \frac{g(x,2)}{f^{m+2}}$. \qed

\begin{rem}
The theorem and its proof generalize nearly 
 verbatim to the situation where $M$ is any
finitely generated $D$-module. 
\end{rem}
Now let us consider the situation in which $M=D/J$ is holonomic on $U$.
Then $D_v/(\{\del_v\}\cdot D_v +D_v\cdot\{1-fv,\phi(J)\})$ is a finite
$D$-module as it is isomorphic to the module
$M\otimes_RR_f$ which is holonomic by theorem 1.3 of \cite{K2}.  
We can find a generator by computing the $b$-function for
integration along $\del_v$ for $D_v/D_v\cdot\{1-vf,\phi(J)\}$ as in step 2 of algorithm
\ref{loc-alg}. If $k$ is the  largest root  
of the $b$-function
then  $v^k$ is a generator for
$D_v/(\{\del_v\}\cdot D_v +D_v\cdot\{1-fv,\phi(J)\})=M\otimes R_f$ (compare \cite{O-T2},
algorithm 5.4). Thus in order to represent the localization $M\otimes
_RR_f$ as a quotient of $D$ all one needs to do is to find the
annihilator of $v^k$ over $D$. This shows the correctness of our algorithm. 

\begin{rem}
Again, the algorithm generalizes to the non-cyclic situation. Let
$M=D^m/J$, $D^m=\oplus_1^mDe_j$. The
modifications are as follows. 
Compute $m$ separate $b$-functions $b_j(s)$ to the integration of 
$D_v(e_j+\phi(J))/(D_v\cdot\{\phi(J),(1-vf)e_j\})$ along $\del_v$. 
$M\otimes_R R_f$ is generated by the cosets of $v^{k_j}e_j$ in $N$ where
$k_j$ is the largest integer root of $b_j(s)$.
\end{rem}

\section{Categorical explanation of the  algorithm}
Let us now give a more categorical explanation of the validity of our
algorithm.

Decompose $j : K^n \setminus V(f) \rightarrow K^n$ as 
$j = p\circ\iota\circ\phi$;
$$ 
\begin{array}{ccc}
  W           & \stackrel{\iota}{\longrightarrow} & K^{n+1} \\
 {\Big \uparrow \rlap{$\vcenter{\hbox{$\scriptstyle\phi$}}$}}&  
 & {\Big \downarrow \rlap{$\vcenter{\hbox{$\scriptstyle p$}}$}} \\
 K^n \setminus V(f) & \stackrel{j}{\longrightarrow}& K^n \\
\end{array}
$$
where
$\phi : K^n \setminus V(f) \rightarrow W = \{(v,x) \mid v\cdot
f(x) = 1\}$ is  
defined by $x\mapsto (1/f(x),x)$; $\iota : W \into K^{n+1}$ is 
the closed embedding, and $p : K^{n+1} \rightarrow K^n$ is the
natural projection. 
Then we have
\[
 M\otimes_RR_f = \int_j j^{-1}M = \int_p\int_\iota\int_\phi j^{-1}M
\]
by the chain rule of integration functors 
(see, e.g., \cite[p.251, 6.4 Proposition]{Borel}, \cite[1.5.1]{HT}).
Note that the chain rule holds in the derived category
in general, but, in this case,
$\int_j$ is an exact functor since
${\cal D}_{X \stackrel{j}{\leftarrow}U} = {\cal D}_X[1/f]$
is flat over ${\cal D}_U$ and
$j$ is an affine morphism.
$\int_\iota \int_\phi$ is also an exact functor.
Hence, we have only to compute  $0$-th integrals in each step of the
computation of 
the integral functors.
$\int_\phi j^{-1}M$ is obtained by the coordinate transformation
represented by our ring map $\phi$;
consider the map
$$y_i = x_i,\ (i=1, \ldots, n),\  y_{n+1} = 1/f(x).$$
Then, $\del_{x_i} = \del_{y_i} - \frac{f_i}{f^2} \del_{y_i}
= \del_{y_i} - {f_i} v^2 \del_{y_i} $ modulo $ 1-fv =0$,
which commute with each other and define our ring map $\phi$.
$\int_\iota$ is nothing but Kashiwara equivalence corresponding to
$\phi(D)\into D_v$. 
It follows that $\int_\iota\int_\phi j^{-1}M=D_v/(D_v\cdot\{\phi(J),1-fv\})$,
compare also  
Proposition A.1 (p.\ 596) of \cite{Oa1}.  
Integration under the projection $p$ corresponds then to our last step
accomplished by forming the quotient
$\left( \int_\iota\int_\phi j^{-1}M \right)/
  \{\del_v\}\cdot \left(\int_\iota\int_\phi j^{-1}M\right)$.

Since $j^{-1}M$ is holonomic, so is $\int_\iota\int_\phi j^{-1}M$ and
hence the $b$-function for integration is nonzero, thus guaranteeing
termination of the search for generators of $\int_jj^{-1}M$.

\section{Examples}
\begin{ex}
For our first example we take $n=1, x=x_1, K=\C, J=x\del_x+\lambda,
\lambda\in\C, f(x)=x$. 

In this scenario, we have to compute the integral of the module
$D_v/(1-xv,x(\del_x-v^2\del_v)+\lambda)$ along $\del_v$. 
As $D_v\cdot\{1-xv,x(\del_x-v^2\del_v)+\lambda\}=
D_v\cdot\{\del_xx-\del_vv+1+\lambda,1-xv\}$, 
one
sees that the $b$-function is $s(s+1+\lambda)$. Thus 
the largest integer root is either $0$ or $-1-\lambda$, depending on
whether $\lambda$ is a negative integer or not.
So $M\otimes_RR_x$ is generated by
$\bar 1\otimes x^{-(-1-\lambda)-2}$ in the former and $\bar 1\otimes 
x^{-2}$ in the latter case.

If for example $\lambda=-7$ we compute a $b$-function of $s(s-6)$
indicating that $M\otimes_RR_x$ is generated by $\bar 1\otimes
x^{-8}$. Since in this case $M=R$ generated by $x^7$, we conclude that
$M\otimes _RR_f$ is in fact generated by $x^{-8}\cdot x^7=1/x$, as it should.

If on the other hand $\lambda=1/2$ then $b(s)=s(s+3/2)$ and hence the
largest integer root is 0. Thus, $M$ is generated by the germ of
$x^{-\lambda}$ and $M\otimes _RR_x$ is generated by
$v^0$ corresponding to $1/x^{2+\lambda}$. $M$ is
already isomorphic to $M\otimes_RR_x$.
\end{ex}
\begin{ex}
In this example we consider the left ideal $J$ generated by
\begin{eqnarray*}
\del_x(x^2-y^3),\\
\del_y(x^2-y^3).
\end{eqnarray*}
These are annihilators of the function $1/f$,
where $f=x^2-y^3$,
but they do not generate the annihilating ideal (see, e.g., \cite{Oa2}). 
The left ideal $J$ is not holonomic since the characteristic
variety of $J$ is
$V(x^3-y^2) \cup V(\xi_x,\xi_y)$,
whose first component has dimension $3$ in $\C^4 =\{(x,y; \xi_x,\xi_y)\}$.
As to an algorithmic method to get the characteristic variety,
see \cite{Oa0}.
We give now the output of a
computer session using the computer algebra system Kan/sm1 (\cite{T2})
interspersed with comments. We
remark that in this case 
$M$ restricted to $U$ is an $\O_U$-coherent free module of
rank one where $U=\C^2\setminus \var(f)$.
\begin{verbatim}
 /ff
   [ ((x^2-y^3)*(Dx - v^2*2*x*Dt) + 2*x)
     ((x^2-y^3)*(Dy + v^2*3*y^2*Dt) -3*y^2)
     (v*(x^2-y^3)-1)
   ]
 def
\end{verbatim}
(this is $\phi(J)$)
\begin{verbatim}

 sm1>ff [(v)] intbfm ::
 [    $216*s^4+1296*s^3+2586*s^2+1716*s$ ] 
 [[6,1],[s+2,1],[6*s+11,1],[6*s+13,1],[s,1]]
\end{verbatim}
(these are the factors of the $b$-function)
\begin{verbatim}
 sm1>ff [(v)] -2 0 1 intall_s ;
 Completed.
\end{verbatim}
(computing the integration)
\begin{verbatim}
  0-th cohomology:  [    1 , 
    [    -3*x*Dx-2*y*Dy-18 , 
         3*y^2*Dx+2*x*Dy , 
         -2*y^3*Dy+2*x^2*Dy-18*y^2 ]  ] 
 -1-th cohomology:  [    0 , [   ]  ] 
\end{verbatim}
The integration $\int_jj^{-1}M$ is not $\O_X$-coherent, although of
course it is still coherent over the sheaf of differential operators
$\D_X$ on $X=\C^2$.
Since localization is an exact functor, the first cohomology group
(corresponding to the first higher order integration) was known to be
zero.
The 0-th cohomology above coincides with the annihilating ideal of
the function $f^{-3}$. 
\end{ex}
\begin{ex}
Let $n=3$, and 
consider the ideal $J$ generated by the system
\begin{eqnarray*}
(x^3-y^2z^2)^2\del_x&+&3x^2,\\
(x^3-y^2z^2)^2\del_y&-&2yz^2,\\
(x^3-y^2z^2)^2\del_z&-&2y^2z.
\end{eqnarray*}
These operators are annihilators of the exponential function
$e^{1/f}$ where $f(x,y,z)=x^3-y^2z^2$. 
The characteristic variety of $M=D/J$
has six components, defined by the prime ideals $(y,x)$, $(z,x)$,
$(\xi_x,\xi_y,\xi_z)$, $(\xi_y,z,x)$, $(\xi_z,y,x)$ and the ideal
generated by
\begin{eqnarray*}
y\xi_y-z\xi_z, & 2x\xi_x+3z\xi_z , & 8z\xi_x^3+27\xi_y^2\xi_z,\\
8y\xi_x^3+27\xi_y\xi_z^2 ,& 
 -4z^2\xi_x^2+9x\xi_y^2, & -4yz\xi_x^2+9x\xi_y\xi_z ,\\
-4y^2\xi_x^2+9x\xi_z^2, & 2z^3\xi_x\xi_z+3x^2\xi_y^2 ,&
2yz^2\xi_x+3x^2\xi_y,\\
2y^2z\xi_x +3x^2\xi_z ,& 
-yz^3\xi_z+x^3\xi_y,& -y^2z^2+x^3,\\
 &  x^3\xi_y^2-z^4\xi_z^2.
\end{eqnarray*} 
All but
the first two are of dimension three. This implies that the
non-holonomic locus of $M$ is contained in the hypersurface 
$x=0$. 
Here, we used \cite{Asir} to obtain the primary ideal decomposition.

Hence we may apply our algorithm to compute $M\otimes_RR_x$. We remark
that contrary to the previous example in this case $j^{-1}M$ is
holonomic but not
coherent as ${\cal O}_U$-module on $U=X\setminus \var(x)$. Using Kan/sm1
again one obtains the eight operators
\begin{eqnarray*}
&&  -3 y \del_y+3 z \del_z, \  
  -2 x y z^2 \del_x-3 x^3 \del_y-4 y z^2, \  
  -2 x y^2 z \del_x-3 x^3 \del_z-4 y^2 z, \  \\
&&
  6 x z^3 \del_x \del_z+9 x^3 \del_y^2+6 x z^2 \del_x+6 y z^2 \del_y+6 z^3 \del_z
    +12 z^2, \  \\
&&  -6 y^2 z^3 \del_z+4 x^4 \del_x+12 x^3 z \del_z+8 x^3+12, \   \\
&&  6 y z^4 \del_z^2-4 x^4 \del_x \del_y-12 x^3 z \del_y \del_z+18 y z^3 \del_z
    -8 x^3 \del_y-12 \del_y, \  \\
&&  8 x^5 \del_x^2+24 x^4 z \del_x \del_z+18 x^3 z^2 \del_z^2+64 x^4 \del_x
    +102 x^3 z \del_z+80 x^3+24 x \del_x+48, \   \\
&&  -6 z^5 \del_z^3+4 x^4 \del_x \del_y^2+12 x^3 z \del_y^2 \del_z-36 z^4 \del_z^2
    +8 x^3 \del_y^2-36 z^3 \del_z+12 \del_y^2
\end{eqnarray*}
which annihilate the function $x^{-2} e^{1/f}$.
The characteristic variety of this holonomic 
left ideal of $D$ has the same last
four components as the characteristic variety of $M$ while the first
two components (of dimension 4) are replaced by $(\xi_z , \xi_y , x)$ and 
  $(z , y , x)$ (of dimension 3).

\end{ex}

\bigbreak
{\em Acknowledgement}:
This paper was written while the authors attended the special program
``Symbolic computations in geometry and analysis'' at the Mathematical Sciences
Research Institute. 
The authors deeply appreciate the efforts of the organizers for
creating  an opportunity 
of studying together intensively.

\end{document}